\documentclass[a4paper]{amsart}

\numberwithin{equation}{section}

\theoremstyle{remark}

\title[Subclasses of Meromorphic Starlike Functions]{Subclasses of Meromorphic Starlike Functions
 Connected to Multiplier Family}
\author[Aghalary, Ebadian and   Eshaghi Gordji ]{R. Aghalary$^1$, A. Ebadian$^2$ and M. Eshaghi Gordji$^3$ }
\address{1,2: Department  of Mathematics, Urmia University,  Urmia, Iran}
 \email{raghalary@yahoo.com and a.ebadian@mail.urmia.ac.ir}
\address{3: Department of Mathematics, Semnan University , P. O. Box 35195-363,
Semnan, Iran}

 \email{madjid.eshaghi@gmail.com}



\begin{document}
\begin{abstract}  The object of this
paper is studying some properties of meromorphic functions which
satisfy in the condition \[Re(zf(z)) > \alpha|z^2f'(z)+zf(z)| .\]
Parallel results for some related classes are also obtained.
\end{abstract}

\maketitle

\section{Introduction and definitions }Denote by $\sum $  the family of functions
$$f(z)={z}^{-1}+\sum_{n=0}^{\infty} a_{n}z^{n} \eqno (1)$$ which are analytic in the punctured disc $ E=\{z : 0<|z| < 1 \}$ with simple pole at $z=0$.
A function $f \in \sum $ is said to be in the class $\sum^{*}
(\alpha)$ of meromorphic starlike functions of order $\alpha$ if
and only if $$ Re \left(\frac{zf'(z)}{f(z)} \right )< -\alpha
\qquad \qquad (z\in E ; 0\leq \alpha < 1). $$ Set $\sum^{*} (0)
=\sum^{*}$.

\vskip .1 truein

We further let $ME(\alpha), 0 \leq \alpha $, be the subclasses of
$\sum$ consisting of functions of the form (1) which satisfy the
condition $$Re(zf(z)) > \alpha|z^2f'(z)+zf(z)|.$$ Set ME=ME(1).
Also for, $0 \leq \alpha <1$, let $MF(\alpha) $ be the subclasses
of $\sum$ consisting of functions of the form (1) which satisfy
the condition

$$ |\frac {zf'(z)}{f(z)} +1| < 1-\alpha , $$

\vskip .1 truein

for $z \in E $. Set $MF=MF(0).$

\vskip .1 truein

Many important properties and characteristics of various
interesting subclasses of the class $\sum$, including (for
example) the class $\sum^{*} (\alpha)$, were investigated by
(among others ) Liu and Srivastava [5],[6] and Aouf [1]. Also some
interesting properties of analytic functions related to multiplier
family were studied by Fournier et al in [3] and Ahuja et al in
[2] and Rosy et al in [9]. In this paper we aim to obtain several
properties of functions belong to the classes $ME(\alpha)$,
$MF(\alpha)$ and $MTE(\alpha)$.

\section{Main Results}
We begin by proving inclusion relation between classes which are
defined in the above.

\vskip .1 truein {\bf Theorem 2.1 }$$ ME(\alpha) \subset
MF(1-\frac{1}{\alpha}) \subset \sum^{*}(1-\frac{1}{\alpha}) \qquad
\qquad 1 \leq \alpha .$$

If $\alpha=1$ all inclusions are proper and for $\alpha >1$ the
result is sharp.

{\bf Proof}. If $f\in ME(\alpha)$, then $$|zf(z)| >
\alpha|z^{2}f'(z)+zf(z)|\quad or \quad |\frac{zf'(z)}{f(z)}+1| <
\frac{1}{\alpha}. \eqno (2)$$

Hence $$Re \frac{zf'(z)}{f(z)} < -1+ \frac{1}{\alpha} \quad
or\quad -Re \frac{zf'(z)}{f(z)} > 1-\frac{1}{\alpha}.\eqno(3)$$ By
making use of (2) and (3) we get our result. But for $\alpha=1$ it
is easy to see that $\frac{e^{z}}{z} \in MF-ME$ and
$\frac{(1-z)^{2}}{z} \in \sum^{*}-MF.$ Now for sharpness set
$f(z)=\frac{1+cz}{z-cz^2}, c=(1+\alpha^2)^{\frac{1}{2}}-\alpha.$
Then $f\in ME(\alpha)$ because for $|z|=r <1$
$$Rezf(z)=\frac{1-c^2r^2}{|1-cz|^2} \geq
\alpha\frac{2cr}{|1-cz|^2}=\alpha|z^2 f'(z)+zf(z)|.$$Note that
$$\frac{-zf'(z)}{f(z)}=1-\frac{2c}{1-c^2z^2},$$ for $z=-r , r
\mapsto 1$, this last expression approaches to
$1-\frac{2c}{1-c^2}=1-\frac{1}{\alpha}.$ Thus $ f \notin
MF(\beta)$ and $f \notin \sum^{*} (\beta)$ for $\beta > 1- \frac
{1}{\alpha} .$

\vskip .1 truein

Next we determine a sufficient condition for a function of the
form (1) to be in the class $ME(\alpha)$.

{\bf Theorem 2.2}. A sufficient condition for a function of the
form (1) to be in the $ME(\alpha)$ is that

$$\sum_{n=0}^{\infty}\left[1+\alpha(n+1) \right] |a_n| \leq 1 .$$

{\bf Proof}. Let $f(z)=z^{-1}+\sum_{n=1}^{\infty}a_nz^n$, then
$$Re(zf(z))=Re \left (1+\sum_{n=0}^{\infty}a_nz^{n+1} \right) \geq
1-\sum_{n=0}^{\infty}|a_n| \eqno(4)$$ and also
$$\alpha|z^2f'(z)+zf(z)| \leq \alpha
|\sum_{n=0}^{\infty}(n+1)a_nz^{n+1}|. \eqno (5)$$By making use of
(4) and (5) we get our result.

\vskip  .1  truein

{Remark 1}. By Theorem 2.1 it follows that $ME \subset \sum^{*}$,
also we note that Theorem 2.2 implies that
$g(z)=z^{-1}+\frac{1}{n+2}z^n \in ME$ for any $n \geq 1$, but if
$n > [\frac {2-3\alpha}{\alpha}]$ then $g$ is not in
$\sum^{*}(\alpha)$, hence $ ME$ is not subset of $
\sum^{*}(\alpha)$ for any $\alpha
> 0$.

\vskip .1 truein

We shall need the following lemma, which is due to Miller and
Mocanu [7] to prove the coefficient estimates for functions
belonging to the class $ME(\alpha)$.

\vskip .1 truein

{\bf Lemma 1}. Let a function $w(z) =a+w_m z^m+\cdots$ be analytic
in the unit disc with $w(z) \neq a$ and $m\geq 1$. If
$z_0=r_0e^{i\theta}$ $(0<r_0<1)$ and $ |w(z_0)|=\max_{|z|\leq
r_0}|w(z)|$. Then $z_0w'(z_0)=kw(z_0) $   and $\Re
\left(1+\frac{z_0 w''(z_0)}{w'(z_0)}\right) \geq k $, where $k$ is
real and $k\geq m$.

\vskip  .1  truein

{\bf Theorem 2.3}. If the function $f$ given by (1) belongs to the
class $ME(\alpha)$, then $$|a_n| \leq
\frac{2}{\sqrt{\alpha^{2}(n+1)^{2}+1}+\alpha(n+1)}, \qquad n \geq
0 .\eqno (6)
$$

The result is sharp for function
$zf(z)=\frac{1+d_{n}z^n}{1-d_{n}z^n}$ where  $d_{n}=
{\sqrt{\alpha^{2}n^{2}+1}+\alpha n}.$

{\bf Proof}. Let $f \in ME(\alpha)$ and $zf(z)=1+Az^{n}+...$. It
is sufficient to show that $|A| \leq 2d_{n}$. For this let
$zf(z)=\frac{1+d_{n}w(z)}{1-d_{n}w(z)}$. It is easy to see that
$w$ is analytic in the unit disc and $w(0)=0$. We wish to show
that $|w(z)| <1 $, for all $z$ in the unit disc. For ,if not, by
Lemma 1 there exists $z_{0}$ in the unit disc such that
$|w(z_{0})| =1$ and $z_{0}w'(z_{0})=kw(z_{0}) , k \geq n$ and
hence $$ Re (z_{0}f(z_{0})) - \alpha
|z_{0}(z_{0}f(z_{0}))'|=\frac{1-d_{n}^2
-2kd_{n}}{|1-d-{n}w(z_{0})|^2} \leq \frac{1-d_{n}^2
-2nd_{n}}{|1-d-{n}w(z_{0})|^2}=0 ,\eqno(7)$$ which contradicts $f
\in ME(\alpha)$. Now the result follows from the well known result
of Robertson [10].

\vskip  .1  truein

{\bf Remark 2}. Also we note that if $f \in ME(\alpha) ,\alpha
\geq 1$,then $Re(z^2f'(z))<0 ,z\in E$. Since if
$f(z)=\frac{g'(z)}{z} \in ME(\alpha)$ where $g(z)$ is an analytic
function in the unit disk, then $f'(z)=\frac{zg''(z)-g'(z)}{z^2}$
. Hence
$$Re(-z^2f'(z))=Reg'(z)-Rezg''(z)>Reg'(z)-\alpha|zg''(z)|>0,$$
which yields result.

\medskip
{\bf 3. Neighborhoods And Partial Sums}

\vskip  .1  truein

Following the earlier works (based upon the familiar concept of
neighborhoods of analytic functions)by Goodman [4] and Ruscheweyh
[8],we begin by introducing here the $\delta-$neighborhood of a
function of the form (1) by the means of the definition

\vskip  .1  truein

 $N_{\delta}(f)=\{g(z)=z^{-1}+\sum_{k=1}^{\infty}b_k
z^k | $ g is analytic in E and $ \sum_{k=1}^{\infty}k|a_k-b_k|\leq
\delta \}.$

\vskip  .1  truein \noindent From Theorem 2.2 for the function
$f(z)=\frac{1}{z}$, we immediately have $N_{\frac
{1}{1+2\alpha}}(f)\subset ME(\alpha).$

\vskip  .1  truein For function $f\in \sum$ given by (1) and $g\in
\sum$ given by $$g(z)=z^{-1}+\sum^{\infty}_{n=0}b_nz^n,$$ we
define the Hadamard product (or convolution) of $f$ and $g$ by
$$ (f*g)(z)=z^{-1}+\sum^{\infty}_{n=0}a_nb_nz^n=(g*f)(z).$$We next give a multiplier convolution characterization for
$ME(\alpha).$

\vskip  .1  truein

{\bf Theorem 3.1} $f\in ME(\alpha) \Leftrightarrow Re
z(f(z)*\frac{1+z(\alpha e^{i \gamma}-1)}{z(1-z)^2 }) > 0 $,
 $\gamma \in(-\pi,\pi], z \in E.$

{\bf Proof}. We have

\begin{eqnarray} {z f(z)+e^{i \gamma}\alpha z(z f)'} &&
= z f(z)*[\frac{1}{1-z} +e^{i\gamma}\alpha\frac{z}{(1-z)^2}] \nonumber \\
& & =z(f(z)*[\frac{1+z(e^{i \gamma}\alpha-1)}{z(1-z)^2}]).
\end{eqnarray}

\noindent Hence we get our result.

{\bf Theorem 3.2}. If $\frac{f(z)-\epsilon z^{-1}}{1-\epsilon} \in
ME(\alpha)$, for $\delta<\epsilon < 1$, then $N_{\gamma}(f)
\subset ME(\alpha)$ where $\gamma=\frac{1}{1+2\alpha} .$

{\bf Proof} . Let $h(z)=z^{-1}+ \sum_{k=0}^{\infty}c_k
z^k=\frac{1+z(\alpha e^{i \gamma}-1)}{z(1-z)^2} .$ It is not
difficult to verify that $|c_k| \leq (1+\alpha(k+1))
,k=0,1,2,3,....$ Let $g \in N_{\gamma}(f)$ and $g(z)=z^{-1}+
\sum_{k=1}^{\infty}b_k z^k .$ Then

$$Re(z(g*h))=Re(z((g-f)*h) +z(f*h))=Re(z((g-f)*h)) +Re(z(f*h)) \eqno  (8)$$

But
$$Re(z((g-f)*h))\geq -|z((g-f)*h)|= -|\sum_{k=0}^{\infty}(b_k-a_k)c_k z^k| > -\delta  ,\eqno(9)$$
since $g \in N_\gamma(f) $. Again $\frac{f- \epsilon z^{-1}}{1-
\epsilon} \in ME(\alpha)$, for $\delta <\epsilon <1$, implies that
$Re (z(\frac{f-\epsilon z^{-1}}{1- \epsilon} *h)) >0 $ by the
Theorem 3.1. That is, $$Re(z(f*h)) >\epsilon  \qquad \qquad  for
\qquad \quad \delta < \epsilon <1 . \eqno(10) $$ Using (9),(10) in
(8) we see that $Re (z(g*h)) >0 $ for all $z \in E.$ Hence Theorem
3.1 show that $g \in ME(\alpha) .$

\medskip
{\bf 4.Negative Coefficients }

In this section at first we introduce the subclass $TME(\alpha)$
consisting of all functions $f\in ME(\alpha)$ which are in the
form
$$f(z) = z^{-1}- \sum_{n=1}^{\infty}a_n z^n
 \qquad \qquad (a_n \geq 0),$$ and then we obtain several
 properties of functions belong to $TME(\alpha)$.

\vskip  .1  truein

{\bf Theorem 4.1}. A function $f$ of the form $f(z) = z^{-1}-
\sum_{n=1}^{\infty}a_n z^n$ is in $TME(\alpha)$ if and only if
$$ \sum_{n=1}^{\infty}(1+\alpha(n+1))a_n \leq 1 .$$ The result
is sharp for the function $f(z)$ given by $$f(z) =
z^{-1}-(\frac{1}{1+\alpha(n+1)})z^n , \qquad \qquad \qquad
 n=1,2,3,....$$

{\bf Proof}. In view of Theorem 2.2 , we need only show that $f
\in TME(\alpha) $ satisfies the coefficient condition . For
$z=re^{i\theta}, 0 \leq r <1 $ and $0 \leq \theta <2\pi$ we have
$rf(r)=1-\sum_{n=1}^{\infty}a_nr^{n+1} $ and $\alpha|r^2f'(r)
+rf(r)|= \alpha \sum_{n=1}^{\infty}(n+1)a_n r^{n+1} .$ The result
follows upon letting $r \mapsto 1 $.

The coefficient characterization of Theorem 4.1 enables us to
determine extreme points and distortion theorems.

\vskip  .1  truein

{\bf Corollary 1}. The extreme points of $TME(\alpha)$ are $
f_{1}(z) = z^{-1} $ and $$f_{n}(z) = z^{-1}
-\frac{z^n}{1+\alpha(n+1)} , n=1,2,3 ,...\quad.$$And $f \in
TME(\alpha)$ if and only if $f$ can be written in the form $$ f(z)
= \sum_{k=1}^{\infty} \lambda_{k}f_{k}(z) ,\qquad where \qquad
\lambda_{k} \geq 0 , \qquad \sum_{k=1}^{\infty} \lambda_{k}=1 .$$

\vskip  .1  truein

{\bf Corollary 2} . If $f(z)=z^{-1}-\sum_{n=1}^{\infty}a_n z^n
,a_{n} \geq 0 $ is in $TME(\alpha),$ then
$$\frac{1}{r}-\frac{1}{1+2\alpha}r \leq |f(z)| \leq
\frac{1}{r}+\frac{1}{1+2\alpha}r ,$$ with equality for
$f(z)=\frac{1}{z}-\frac{1}{1+2\alpha}z $ at $z=r,ir.$

\vskip  .1  truein

Finally we prove

{\bf Theorem 4.2} . Let $f\in \sum$ be given by (1) and define the
partial sums $S_1 (z)$ and $S_n(z)$ by $S_1 (z)=z^{-1}$ and
$S_n(z)= z^{-1} +\sum _{k=1}^{n-1} a_k z^k$.

Suppose also that $$ \sum_{K=1}^{\infty} d_k|a_k| \leq 1   \qquad
\qquad  (d_k=1+\alpha(n+1)). \eqno (11)$$

Then we have $$ Re(\frac{f(z)}{S_n(z)}) > 1-\frac{1}{d_n} \qquad
and \qquad Re(\frac{S_n(z)}{f(z)}) > \frac {d_n}{1+d_{n}} \qquad
(z\in E ; n\in N=\{1,2,3,...\}) .\eqno (12)$$

Each of the bounds in (12) is the best possible for $n\in N$

{\bf Proof} . For the coefficients $d_k$ given by (11), it is not
difficult to verify that $ d_{k+1} >d_{k} >1, k=1,2,3,... .$
Therefore, by using the hypothesis (11), we have $$
\sum_{k=1}^{n-1}|a_k| + d_n\sum_{k=n}^{\infty}|a_k| \leq
\sum_{k=1}^{\infty}d_k|a_k| \leq 1 . \eqno (13)$$ By  setting $$
g_1(z)= d_n ( \frac{f(z)}{s_n(z)}-(1- \frac{1}{d_n}))=1+\frac{d_n
\sum_{k=n}^{\infty} a_k z^{k+1}}{1+\sum_{k=1}^{n-1} a_k z^{k+1}}
\eqno (14)$$ and applying (13), we find that $$ \left|
\frac{g_1(z) -1}{g_1(z)+1} \right | \leq \frac{
d_n\sum_{k=n}^{\infty} |a_k|}{2-2\sum_{k=1}^{n-1}|a_k| -
d_n\sum_{k=n}^{\infty}|a_k|}
 \leq 1  \quad (z \in E),\eqno(15)$$which readily yields the left
 assertion (12) of Theorem 4.2. If we take $$ f(z)= z^{-1}-
 \frac{z^{n}}{d_n}, \eqno (16)$$ then
 $$\frac{f(z)}{S_n(z)}=1-\frac{z^{n+1}}{d_n} \mapsto 1-\frac{1}{d_n} \quad as \quad z \mapsto 1^{-} ,$$ which shows that the bound in (12) is the best possible
 for each $n\in N$. Similarly, if we put $$ g_2(z)= (1+d_n) \left (
 \frac {S_n(z)}{f(z)}-\frac{d_n}{1+d_n} \right )
 =1-\frac{(1+d_n)\sum_{k=n}^{\infty}a_k z^{k+1}}{1+\sum_{k=1}^{\infty}a_k
 z^{k+1}} \eqno  (17) $$ and make use of (13) we obtain $$ \left |
 \frac{g_2(z)-1}{g_2(z)+1} \right| \leq
 \frac{(1+d_n)\sum_{k=n}^{\infty}|a_k|}{2-2\sum_{k=1}^{n-1}|a_k|
  +(1-d_n)\sum_{k=n}^{\infty}|a_k|} \leq 1 \quad (z \in E) $$ which
  leads us to the assertion (12) of Theorem 4.2. The bounds
  given in the right of (12) is sharp with the function given by
  (16). The proof of Theorem 4.2 is thus complete.

\end{document}